\documentclass{birkjour}

\usepackage{amsfonts}
\usepackage{amsmath}
\usepackage{amsthm}
\usepackage{amssymb}
\usepackage[english]{babel}
\usepackage{graphics}
\usepackage{latexsym}
\usepackage{longtable} 
\usepackage{mathrsfs} 
\usepackage{graphicx}
\usepackage{wrapfig} 
\usepackage[pdftex,colorlinks=true,linkcolor=blue]{hyperref}
\usepackage{enumitem}
\usepackage{esint}

\newtheorem{theorem}{Theorem}[section]
\newtheorem{corollary}[theorem]{Corollary}
\newtheorem{lemma}[theorem]{Lemma}
\newtheorem{proposition}[theorem]{Proposition}

\newtheorem{claim}[theorem]{Claim}

\theoremstyle{definition}
\newtheorem{definition}[theorem]{Definition}
\theoremstyle{remark}
\newtheorem{remark}[theorem]{Remark}
\newtheorem{example}{Example}

 \numberwithin{equation}{section}

\newcommand{\mL}{\mathcal{L}}
\newcommand{\mH}{\mathcal{H}}
\newcommand{\mF}{\mathcal{F}}

\newcommand{\mP}{\mathscr{P}}

\newcommand{\mO}{\mathcal{O}}

\newcommand{\mD}{\mathcal{D}}

\newcommand{\E}{\mathrm{E}}
\newcommand{\D}{\mathrm{D}}
\renewcommand{\H}{\mathrm{H}}

\newcommand{\R}{\mathbb{R}}

\newcommand{\N}{\mathbb{N}}

\newcommand{\noi}{\noindent}
\newcommand{\ms}{\medskip}

\newcommand{\de}{\delta}
\newcommand{\De}{\Delta}
\newcommand{\e}{\varepsilon}

\newcommand{\la}{\lambda}
\newcommand{\La}{\Lambda}

\newcommand{\Om}{\Omega}
\newcommand{\om}{\omega}

\newcommand{\av}{-\!\!\!\!\!\!\int}

\newcommand{\weak }{\, -\!\!\!\!-\!\!\!\rightharpoonup}
\newcommand{\weakstar }{ \overset{\, *_{\phantom{|}}}{{\smash{\weak }}\, } }

\newcommand{\larrow}{\longrightarrow}
\newcommand{\ot}{\otimes}

\newcommand{\ri}{\rightarrow}

\newcommand{\p}{\partial}
\newcommand{\sub}{\subseteq}
\newcommand{\set}{\setminus}
\newcommand{\by}{\times}

\newcommand{\ess}{\mathrm{ess}}
\newcommand{\diam}{\mathrm{diam}}
\newcommand{\dist}{\mathrm{dist}}

\newcommand{\bt}{\begin{theorem}}\newcommand{\et}{\end{theorem}}
\newcommand{\bd}{\begin{definition}}\newcommand{\ed}{\end{definition}}
\newcommand{\bl}{\begin{lemma}}\newcommand{\el}{\end{lemma}}
\newcommand{\beq}{\begin{equation}}\newcommand{\eeq}{\end{equation}}
\newcommand{\bc}{\begin{claim}}\newcommand{\ec}{\end{claim}}
\newcommand{\bex}{\begin{example}}\newcommand{\eex}{\end{example}}
\newcommand{\bcor}{\begin{corollary}}\newcommand{\ecor}{\end{corollary}}
\newcommand{\bp}{\begin{proof}}\newcommand{\ep}{\end{proof}}

\newcommand{\BPL}{\medskip \noindent \textbf{Proof of Lemma} }

\newcommand{\BPCOR}{\medskip \noindent \textbf{Proof of Corollary} }
\newcommand{\BPP}{\medskip \noindent \textbf{Proof of Proposition} }
\newcommand{\BPT}{\medskip \noindent \textbf{Proof of Theorem} }

\numberwithin{equation}{section}

\begin{document}

%
%
%
%
%
%
%
%
%

\title[The Eigenvalue Problem for the $\infty$-Bilaplacian]{The Eigenvalue Problem for the $\infty$-Bilaplacian}

\author{Nikos Katzourakis}

\address{Department of Mathematics and Statistics, University of Reading, Whiteknights, PO Box 220, Reading RG6 6AX, UNITED KINGDOM}
\email{n.katzourakis@reading.ac.uk}

\thanks{N.K. has been partially financially supported through the EPSRC grant EP/N017412/1}
  
\author{Enea Parini}

\address{Aix Marseille Univ, CNRS, Centrale Marseille, I2M, 39 Rue Frederic Joliot Curie, 13453 Marseille, FRANCE}
\email{enea.parini@univ-amu.fr}



\date{June 15, 2017}

\begin{abstract} We consider the problem of finding and describing minimisers of the Rayleigh quotient
\[
\La_\infty \, :=\, \inf_{u\in \mathcal{W}^{2,\infty}(\Om)\set\{0\} }\frac{\|\De u\|_{L^\infty(\Om)}}{\|u\|_{L^\infty(\Om)}},
\]
where $\Omega \subseteq \mathbb{R}^n$ is a bounded $C^{1,1}$ domain and $\mathcal{W}^{2,\infty}(\Om)$ is a class of weakly twice differentiable functions satisfying either $u=0$ on $\partial \Om$ or $u=|\D u|=0$ on $\partial \Om$. Our first main result, obtained through approximation by $L^p$-problems as $p\to \infty$, is the existence of a minimiser $u_\infty \in \mathcal{W}^{2,\infty}(\Om)$ satisfying
\[
\left\{
\begin{array}{ll}
\De u_\infty \, \in \, \La_\infty \mathrm{Sgn}(f_\infty) & \text{ a.e.\ in }\Om,
\ms
\\
\De f_\infty \, =\, \mu_\infty &  \text{ in }\mD'(\Om), 
\end{array}
\right.
\]
for some $f_\infty\in L^1(\Om)\cap BV_{\text{loc}}(\Om)$ and a measure $\mu_\infty \in \mathcal{M}(\Om)$, for either choice of boundary conditions. Here Sgn is the multi-valued sign function. We also study the dependence of the eigenvalue $\La_\infty$ on the domain, establishing the validity of a Faber-Krahn type inequality: among all $C^{1,1}$ domains with fixed measure, the ball is a strict minimiser of $\Om \mapsto \La_\infty(\Om)$. This result is shown to hold true for either choice of boundary conditions and in every dimension.

\end{abstract}

\maketitle


\section{Introduction} \label{section1}

Let $\Omega \subseteq \mathbb{R}^n$ be a bounded domain with $C^{1,1}$ boundary $\p\Om$, where $n\in \N$. In the present paper we are interested in studying nonlinear higher order $L^\infty$ eigenvalue problems. More precisely, we consider the problem of existence of minimisers to the $L^\infty$ Rayleigh quotient
\[
\La_\infty(\Om)\, :=\, \inf_{u\neq 0 }\frac{\|\De u\|_{L^\infty(\Om)}}{\|u\|_{L^\infty(\Om)}}
\]
over appropriate classes of twice weakly differentiable functions, involving two distinct types of boundary conditions on $\p\Om$. Furthermore, we are interested in studying the structure of these minimisers as well as the dependence of the eigenvalue $ \La_\infty(\Om)$ on the shape of the domain $\Om$. 

The types of boundary conditions we will consider are either Dirichlet conditions ($u=0$ on $\p\Om$), which we refer to as the ``hinged case", or coupled Dirichlet-Neumann conditions ($u=|\D u|=0$ on $\p\Om$) which we will refer to as the ``clamped case". The respective hinged and clamped functional spaces wherein we will minimise the $L^\infty$ Rayleigh quotient are
\begin{align}
\mathcal{W}^{2,\infty}_{\mathrm{H}}(\Omega) \, &:= \, \bigcap_{1<p<\infty} \Big\{ u \in \big(W^{2,p} \cap W^{1,p}_0\big)(\Om)\ :\ \De u \in L^\infty(\Om)\Big\},
\label{1.1}
\\
\mathcal{W}^{2,\infty}_{\mathrm{C}}(\Omega) \, &:= \, \bigcap_{1<p<\infty} \Big\{ u \in W^{2,p}_0(\Om)\ :\ \De u \in L^\infty(\Om)\Big\},
\label{1.2}
\end{align}
and they are Fr\'echet spaces. Our general notation will be either a convex combination of standard symbolisations (as e.g.\ in \cite{E,EG,B}) or else self-explanatory and clear from the context. The respective hinged and clamped eigenvalues will be denoted as
\begin{align}
\label{1.3} 
\Lambda_\infty^{\mathrm{H}}(\Omega)\,  &:=  \inf_{u \in \mathcal{W}^{2,\infty}_{\mathrm{H}}(\Omega)\setminus \{0\}} \frac{\|\Delta u\|_{L^\infty(\Omega)}}{\|u\|_{L^\infty(\Omega)}},
\\
\Lambda_\infty^{\mathrm{C}}(\Omega) \, & := \inf_{u \in \mathcal{W}^{2,\infty}_{\mathrm{C}}(\Omega)\setminus \{0\}} \frac{\|\Delta u\|_{L^\infty(\Omega)}}{\|u\|_{L^\infty(\Omega)}},
\label{1.4}
\end{align}
and the dependence of $\Lambda_\infty^{\mathrm{H}}(\Om),\Lambda_\infty^{\mathrm{C}}(\Om)$ on the domain $\Om$ will be suppressed if it is fixed and we do not vary it. This problem can be seen as the higher-order generalisation of the eigenvalue problem for the $\infty$-Laplacian, which has been first studied by Juutinen, Lindqvist and Manfredi in \cite{JLM}. Inspired by their results, herein we prove existence of $\infty$-eigenfunctions by approximation, considering the respective $L^p$-Rayleigh quotients
\begin{align} 
\Lambda_p^{\mathrm{H}}(\Omega)\,  &:=  \inf_{u \in ({W}^{2,p}\cap W^{1,p}_0)(\Omega)\setminus \{0\}} \frac{\|\Delta u\|_{L^p(\Omega)}}{\|u\|_{L^p(\Omega)}},
\label{1.5}
\\
\Lambda_p^{\mathrm{C}}(\Omega)\,  &:=  \inf_{u \in {W}^{2,p}_0(\Omega)\setminus \{0\}} \frac{\|\Delta u\|_{L^p(\Omega)}}{\|u\|_{L^p(\Omega)}},
\label{1.6}
\end{align}
for $p\in(1,\infty)$ and taking $p\to \infty$. By standard weak compactness, lower semicontinuity and Lagrange multiplier arguments, one easily sees that (for finite $p$) minimisers $u_p$ to the respective $L^p$-eigenvalue problems do exist and solve the Dirichlet problems
\begin{equation}
\label{hinged} 
\left\{
\begin{array}{l l l} \Delta^2_p\, u_p \,=\, (\Lambda_p)^p |u_p|^{p-2}u_p & \text{in }\Omega,
\\ 
\ \ \ \ \, u_p \,=\, 0  & \text{on }\partial \Omega,
\\
\ \  \De u_p\, =\, 0 & \text{on }\partial \Omega,
\end{array}
\right. 
\end{equation}
where $(u_p,\La_p)=(u_p^{\mathrm{H}}, \Lambda_p^{\mathrm{H}})$ and
\begin{equation}
\label{clamped} 
\left\{
\begin{array}{l l l} \Delta^2_p \,u_p \,=\, (\Lambda_p)^p |u_p|^{p-2}u_p & \text{in }\Omega,
\\ 
\ \ \ \ \, u_p \,=\, 0 & \text{on }\partial \Omega,
\\
\ \  \D u_p\, =\, 0 & \text{on }\partial \Omega,
\end{array}
\right.
\end{equation}
where $(u_p,\La_p)=(u_p^{\mathrm{C}}, \Lambda_p^{\mathrm{C}})$. In the above $\De^2_p$ is the $p$-Bilaplace operator, given by
\[ 
\Delta^2_p \,u \,:= \, \Delta (|\Delta u|^{p-2}\Delta u).
\]
The eigenvalue problem for the $p$-Bilaplacian, apart from the linear case $p=2$, has not received much attention thus far. In the linear case, the hinged eigenvalue problem is not very meaningful, because the first eigenvalue is simply given by the square of the first eigenvalue of the Laplacian under Dirichlet boundary conditions. For $p \neq 2$, Dr\'abek and Otani showed in \cite{drabekotani} that the first eigenfunction is unique (up to a multiplicative constant) and strictly positive (or negative) inside $\Omega$. Furthermore, as a straightforward application of Talenti's symmetrisation principle \cite{talenti}, which we recall in our second Appendix, a Faber-Krahn type inequality holds true: among all domains with fixed volume, the first eigenvalue is minimised by the ball up to perhaps rigid motions.

On the other hand, the clamped eigenvalue problem presents several interesting features already in the case of $p=2$, which make its study a highly nontrivial matter. Indeed, the first eigenfunction might be sign-changing, even for relatively simple domains such as squares or elongated ellipses \cite{coffmansquare}. Moreover, some domains admit more than one first eigenfunction, as shown in \cite{coffmanduffin}. However, if $\Omega$ is a ball, the first eigenfunction is unique and strictly positive (see for instance \cite[Theorem 3.7]{gazzolagrunausweers}). The Faber-Krahn inequality has been shown to hold true only in dimensions $n=2$ \cite{nadirashvili} and $n=3$ \cite{ashbaughbenguria}, while it still remains a challenging open problem in higher dimensions.  The limiting case $p=1$ has been studied by the second author jointly with Ruf and Tarsi in \cite{pariniruftarsi, pariniruftarsi2}, wherein results analogous to the case $p=2$ were obtained. However, in the clamped case, positivity of the first eigenfunction in a ball and the Faber-Krahn inequality were shown to be true only in dimension $n=2$.

The first main result in the present work concerns the existence and the qualitative structure properties of minimisers and is given below. 

\smallskip

\noi \textbf{Notational convention:} \emph{For the sake of simplicity and to avoid repetition, we will drop the sub/superscripts $``\,\mathrm{C},\mathrm{H}$" and treat both cases in a unified fashion, indicating any differences between the clamped and the hinged case where appropriate.}

\bt[Existence, structure and approximation for the $\infty-$ei\-gen\-prob\-lem] 
\label{theorem1}
Let $\Om \sub \R^n$ be a bounded domain with $C^{1,1}$ boundary and consider the $L^\infty$ variational problems \eqref{1.3} and \eqref{1.4}, placed in the respective spaces \eqref{1.1} and \eqref{1.2}. For $1<p<\infty$, consider the pair of $L^p$ normalised eigenfunction-eigenvalue, corresponding to either \eqref{hinged}-\eqref{1.5} or \eqref{clamped}-\eqref{1.6}:
\[
(u_p,\La_p) \in \big(W^{2,p}(\Om)\set\{0\}\big) \by (0,\infty), \ \ \ \|u_p\|_{L^p(\Omega)}=1.
\]
Then, there exists a sequence of exponents $(p_\ell)_{\ell =1}^\infty$ tending to infinity, such that
\[
\left\{
\begin{split}
& (u_{p_\ell},\La_{p_\ell}) \larrow (u_\infty,\La_\infty), \hspace{17pt} \text{ in }C^{1}(\overline{\Om})\by \R, 
\\
& \D^2 u_{p_\ell} \weak \D^2 u_\infty,  \hspace{45pt} \text{ in }L^q(\Om,\R^{n\by n}) \text{ for all }q \in (1,\infty),
\ms
\\
& \dfrac{|\De u_{p_\ell}|^{{p_\ell}-2}\De u_{p_\ell}}{(\La_{p_\ell})^{p_\ell}} 
\, \bigg\{
\begin{array}{l l}
\!\!\larrow f_\infty, & \text{ in }L^q_{\text{loc}}(\Om)\text{ for all }q \in \Big[1,\dfrac{n}{n-1}\Big),
\ms
\\
\!\!\weakstar\, f_\infty, & \text{ in }BV_{\text{loc}}(\Om),
\end{array}
\ms
\\
& |u_{p_\ell}|^{{p_\ell}-2}u_{p_\ell} \,\mL^n \,\, \weakstar \,\, \mu_\infty, \hspace{21pt} \text{ in }\mathcal{M}(\Om),
\end{split}
\right.
\]
as $\ell \ri\infty$, where $f_\infty\in L^1(\Om)\cap BV_{\text{loc}}(\Om)$, $\mu_\infty \in \mathcal{M}(\Om)$, and $u_\infty \in \mathcal{W}^{2,\infty}(\Om)$ is a normalised minimiser of \eqref{1.3} or \eqref{1.4} respectively, satisfying $\|u_\infty\|_{L^\infty(\Omega)}=1$ and
\[
\La_\infty\, =\,  \|\Delta u_\infty\|_{L^\infty(\Omega)}.
\]
Moreover, $0<\La_\infty<\infty$, and $f_\infty$, $\mu_\infty$ are such that
\beq 
\label{*}
\left\{
\begin{array}{ll}
\De u_\infty(x) \, \in\, \La_\infty \mathrm{Sgn}\big(f_\infty(x)\big) & \text{ a.e.\ }x\in\Om,
\ms
\\
\ \ \ \ \De f_\infty  =\, \mu_\infty & \text{ in }\mD'(\Om).
\end{array}
\right.
\eeq
In the above,  $\mathrm{Sgn} :\R \larrow 2^{\R}$ is the set-valued sign function given by:
\[
\mathrm{Sgn}(x)\,:=\, 
\left\{
\begin{array}{r l}
\{-1\}, & x<0,
\\ 
\phantom{.}[-1,+1], & x=0,
\\
\{+1\}, & x>0.
\end{array}
\right.
\]
In particular, in the case of hinged boundary conditions, one has $\mu_\infty \equiv 0$, $f_\infty \equiv 1$ and $\De u_\infty \equiv \La_\infty$ on $\Om$.

\et 

The symbolisation ``$\mL^n$" above obviously stands for the Lebesgue measure in $\R^n$. Theorem \ref{theorem1} establishes the existence of second order $\infty$-ei\-gen\-func\-tions which solve the parametric system \eqref{*} consisting of a second order differential inclusion (satisfied in the strong sense a.e.\ in $\Om$), coupled by a second order equation with measure right hand side (satisfied in the distributional sense). The system \eqref{*} can be seen as a kind of ``constrained Euler-Lagrange equations" for the $L^\infty$ second order eigenvalue problems \eqref{1.3}-\eqref{1.4}, but we temporarily defer the relation of the present variational problems to the existing theory of Calculus of Variations in $L^\infty$ until later in the introduction, after we will have expounded on our second main result.

In addition to existence and structure, Theorem \ref{theorem1} provides extra information for our $L^\infty$-eigenproblems, showing they are approximable by more conventional $L^p$-eigenvalue problems. In the case of hinged boundary conditions, the weaker requirements on the boundary data allow to have $\infty$-eigenfunctions with constant Laplacian throughout the domain. However, in the clamped case the $\infty$-eigenfunctions are non-$C^2$ even in one space  dimension (see Remark \ref{remark10}). Moreover, it appears that one can not in general expect the differential inclusion to reduce to an equation because the level set $\{f_\infty=0\}$ might have positive measure, as shown by the example $f(x)=\max\{x,0\}$ which solves $f''=\de_0$ in $\mD'(\R)$.

Our second main result concerns an inequality of Faber-Krahn type. Namely, we study the dependence of the eigenvalues \eqref{1.3}-\eqref{1.4} on the geometry of the domain $\Om$, under a volume constraint. The relevant theorem below establishes that the Euclidean ball is a strict minimiser of both $\Om \mapsto \La^{\mathrm{H}}_\infty(\Om)$ and of $\Om \mapsto \La^{\mathrm{C}}_\infty(\Om)$, among all regular bounded domains with fixed measure. Let us stress that our result holds true \emph{in every dimension}, even in the clamped case.

\begin{theorem}[Faber-Krahn inequality for the second order $\infty$-eigenproblem]
\label{theorem2}
Let $\mathrm{B}$ be an open ball in $\R^n$ and let $\Omega$ be a bounded domain in $\R^n$ with $C^{1,1}$ boundary of the same measure as $\mathrm{B}$. Then, in the case of either hinged or clamped boundary conditions, the eigenvalues \eqref{1.3}-\eqref{1.4} satisfy
\begin{equation} \label{faberkrahn}
\ \ \Lambda_\infty (\Omega) \,\geq \, \Lambda_\infty (\mathrm{B}).
\end{equation}
Equality in \eqref{faberkrahn} holds true if and only if $\Omega =\mathrm{B}$, up to rigid motions in $\R^n$.
\end{theorem}

Finally, we study the concrete case when $\Om$ is a ball in $\R^n$ and we calculate explicitly the eigenvalues $\La^{\mathrm{H}}_\infty, \La^{\mathrm{C}}_\infty$ and the eigenfunctions $u^{\mathrm{H}}_\infty, u^{\mathrm{C}}_\infty$ in this case (see Corollary \ref{corollary4}, Proposition \ref{proposition8}). In both cases, the eigenfunctions are unique (up to a multiplicative constant), strictly positive and radially decreasing.

We conclude this introduction by placing the second order $\infty$-eigenvalue problem we study herein in the wider context of Calculus of Variations in $L^\infty$. Variational problems for first order functionals 
\beq  \label{1.11}
\ \ \ \ \   \E_\infty(u,\mO) \,=\,\underset{x\in \mathcal{O}}{\mathrm{ess}\,\sup}\, \H \big(x, u(x),\mathrm{D} u(x)\big) ,\ \ \ u\in W^{1,\infty}(\Om),\ \mO \in \mL(\Omega),
\eeq
together with the associated equations, first emerged in the work of Aronsson in the 1960s (\cite{A1}--\cite{A3}). The area is now well developed and the relevant bibliography is vast; for a pedagogical introduction accessible to non-experts, we refer to \cite{K1} (see also \cite{C}). Higher order $L^\infty$ variational problems have only very recently begun to be investigated and are still poorly understood. In the recent paper \cite{KP2}, the first author jointly with Pryer considered second order variational problems and their relevant equations, focusing on functional of the form
\[
\ \ \ \ \   \E_\infty(u,\mO) \,=\, \underset{x\in \mathcal{O}}{\mathrm{ess}\,\sup}\, \H \big(\mathrm{D}^2 u(x)\big) ,\ \ \ u\in W^{2,\infty}(\Om),\ \mO \in \mL(\Omega).
\]
Subsequently, in a joint paper with Moser \cite{KMo} the case of dependence on second derivatives through the Laplacian was considered, focusing on the model case of so-called \emph{$\infty$-Bilaplacian}: 
\beq \label{1.12}
\De^2_\infty u\, :=\, \big( \De u\, \mathrm{I} \big)^{\ot 3}:  \big(\D^3u \big)^{\ot 2}\, =\, 0.
\eeq

In the light of the above general $L^\infty$ framework, we see the quantities $\Lambda_\infty^{\mathrm{H}}(\Omega)$ and $\Lambda_\infty^{\mathrm{C}}(\Omega)$ as the \emph{first eigenvalues} of the $\infty$-Bilaplacian under the respective (hinged or clamped) boundary conditions and the parametric system \eqref{*} as the analogue of the constrained Euler-Lagrange equations for the minimisation problems \eqref{1.3}-\eqref{1.4}. However, there does exist a more conventional PDE arising in the formal limit of the Dirichlet problems \eqref{hinged}-\eqref{clamped} as $p\to \infty$: by exploiting the relation
\[ 
\Delta^2_p u \,= \, (p-1)|\Delta u|^{p-2}\Delta^2 u \,+\, (p-1)(p-2)|\Delta u|^{p-4}\Delta u |\D(\Delta u)|^2
\]
and performing similar computations as in \cite{JLM}, one can see that any putative $\infty$-eigenfunction $u_\infty$ has to satisfy
\[
\min \big\{|\Delta u|-\Lambda_\infty |u|\, , \,\Delta^2_\infty u \big\} \, =\, 0,
\]
where $\Delta^2_\infty$ is the $\infty$-Bilaplacian given by \eqref{1.12}. Notwithstanding, this is merely a formal claim, since we can not expect the solutions to be classical, and, to the best of our knowledge, there does not exist any analogue of the theory of viscosity solutions for the higher order problem at hand which is equally stable under limiting processes. However, this is not an issue because for the particular problem herein, the method of $L^p$-approximations constructs second order $\infty$-eigenfunctions with finer structure. This renders the direct study of the formal third order PDE redundant, whilst we obtain also a selection principle of the numerous possible $\infty$-eigenfunctions realising the infima in \eqref{1.3}-\eqref{1.4}. A similar phenomenon has already arisen in the paper \cite{KMo}, wherein the authors proved existence and uniqueness of (absolute) minimisers to $u\mapsto \|\De u\|_{L^\infty(\Om)}$ by solving the parametric system
\beq 
\label{**}
\left\{
\begin{array}{ll}
\De u_* \, =\, \La_* \mathrm{sgn}(f_*) & \text{ a.e.\ in }\Om,
\smallskip
\\
\De f_* \, =\, 0 & \text{ a.e.\ in }\Om,
\end{array}
\right.
\eeq
for any given prescribed boundary values $u_*=g$ and $\D u_*=\D g$ on $\p\Om$. In \eqref{**}, ``sgn" is the usual single-valued sign function. In particular, \eqref{**} implies that $|\De u_*|=\La_*$ a.e.\ in $\Om$ and any such $u_*$ is the unique minimising $\infty$-Biharmonic function solving \eqref{1.12} in the appropriate sense of \emph{$\mD$-solutions}, a new theory of generalised solutions for fully nonlinear systems recently introduced in \cite{K2,K3}. The fact that $u_*$ solves \eqref{1.12} if it solves \eqref{**} can be readily seen formally by recasting \eqref{1.12} as $\De u\, |\D \big(|\De u|^2)|^2= 0$.
 \ms

\section{Existence, structure and $p$-approximation to the eigenproblem for the $\infty$-Bilaplacian}

Let $\Om\sub\R^n$ be a given domain with $C^{1,1}$ boundary $\p\Om$. In this section we establish Theorem \ref{theorem1}. Its proof consists of several lemmas and, as in the statement, we tackle both cases simultaneously. To this end, it suffices to consider only the case of hinged boundary conditions, because if we obtain the desired existence-compactness-approximation conclusion by requiring the weaker condition ``$u=0$ on $\p\Om$" for the $L^p$ approximating sequences of eigenfunctions, then it most certainly holds under the stronger requirement ``$u=|\D u|=0$ on $\p\Om$" of clamped boundary conditions. Also, the putative limit eigenfunction $u_\infty$ will be in the respective space because
\[
\mathcal{W}^{2,\infty}(\Om)\,\supseteq \, \mathcal{W}^{2,\infty}_{\mathrm{H}}(\Om)\,\supseteq \, \mathcal{W}^{2,\infty}_{\mathrm{C}}(\Om)
\]
and the hinged/clamped functional spaces are closed in their super-space
\[
\mathcal{W}^{2,\infty}(\Om)\,:=\, \underset{1<p<\infty}{\bigcap}\Big\{ u\in W^{2,p}(\Om) \ : \De u \in L^\infty(\Om) \Big\}.
\]
For technical convenience in the proof we modify our notation slightly, as follows: for $p\in[1,\infty]$, we consider the normalised $L^p$ norm with respect to the probability measure $\la=\mL^n / \mL^n(\Om) \in \mP(\Om)$, that is
\beq
\label{2.1}
\|f\|_{L^p(\Om,\la)}\,:=\, \left\{
\begin{array}{l l}
\displaystyle \left(\av_\Om |f|^p \right)^{1/p}, & 1\leq p <\infty, 
\ms
\\
\|f\|_{L^\infty(\Om)}\, , & p=\infty,
\end{array}
\right.
\eeq
and, given a fixed $p\in(1,\infty)$, we also consider the constrained variational problem of finding $u_p \in \big({W}^{2,p}\cap W^{1,p}_0\big)(\Omega)$ with $\|u_p\|_{L^p(\Omega,\la)}=1$ such that
\beq
\label{2.2}
\|\Delta u_p\|_{L^p(\Omega,\la)} \, =\, \Lambda_p,
\eeq
where
\beq
\label{2.3} 
\Lambda_p\,  := \, \inf \Big\{\|\Delta v\|_{L^p(\Om,\la)} \ : \ v \in \big({W}^{2,p}\cap W^{1,p}_0\big)(\Omega),\, \|v\|_{L^p(\Om,\la)} =1 \Big\}.
\eeq
By standard weak compactness, lower semicontinuity and Lagrange multiplier arguments (see e.g.\ the relevant arguments for the Laplacian in \cite{E}), one easily sees that for any $p \in (1,\infty)$ there indeed exists a desired minimiser $u_p$ of \eqref{2.2}-\eqref{2.3} which solves weakly the Dirichlet problem
\begin{equation}
\label{2.4}
\left\{
\begin{array}{l l l} \Delta^2_p\, u_p \, =\,  (\Lambda_p)^p |u_p|^{p-2}u_p & \text{ in }\Omega,
\ms
\\ 
\ \ \ \ \, u_p \, = \, 0  & \text{ on }\partial \Omega.
\end{array}
\right. 
\end{equation}
Note that we refrain from stating the natural boundary condition ``$\De u_p =0$ on $\p\Om$" which is also satisfied weakly in the hinged case only, because we do not utilise it in any way in the foregoing reasoning which applies to both cases.

\smallskip

We begin with the next lemma.

\bl \label{lemma2.2}
Let $\{(u_p,\La_p) : 1<p<\infty\}$ be the family of eigenfunctions and eigenvalues solving for each $p$ the problems \eqref{2.2}-\eqref{2.4} and such that $\|u_p\|_{L^p(\Omega,\lambda)}=1$. Then, for any sequence of indices $p$ tending to infinity, there exists a subsequence $(p_\ell)_{\ell =1}^\infty$ and a pair
\[
(u_\infty,\hat{\La}_\infty) \,\in \, \mathcal{W}^{2,\infty}(\Om) \by [0,\infty)
\]
with 
\[
\|u_\infty\|_{L^\infty(\Om)} \, =\, 1
\]
such that
\[
\left\{
\begin{split}
& (u_{p},\La_{p}) \larrow (u_\infty,\hat{\La}_\infty), \hspace{10pt} \text{ in }C^{1}(\overline{\Om})\by [0,\infty), 
\\
& \ \ \ \D^2 u_{p} \weak \D^2 u_\infty,  \hspace{24pt} \text{ in }L^q(\Om,\R^{n\by n}) \text{ for all }q \in (1,\infty),
\end{split}
\right.
\]
along this subsequence as $p\to \infty$. In particular, $\De u_\infty \in L^\infty(\Om)$ and we also have
\[
\| \De u_\infty \|_{L^\infty(\Om)}\, =\, \hat{\La}_\infty\,\|u_\infty\|_{L^\infty(\Om)}\, =\, \hat{\La}_\infty .
\]
\el

\BPL \ref{lemma2.2}. Consider any increasing sequence of indices $p$ tending to infinity and suppose we are along this sequence. We begin by obtaining an a priori bound for the sequence $(\La_p)_1^\infty$. Fix $\xi \in C^\infty_c(\Om)$, $\xi \not\equiv 0$. Then, there exists a modulus of continuity $\om \in C(0,\infty)$ with $\om(0^+)=0$ and $0<\om<1/2$ such that
\[
\phantom{_\big|}    
\|\xi\|_{L^\infty(\Om,\la)} \,\geq\, \|\xi\|_{L^p(\Om,\la)} \,\geq \big(1-\om(1/p)\big) \, \|\xi\|_{L^\infty(\Om,\la)}.
\]
Note also that \eqref{2.1} implies $\|\cdot\|_{L^\infty(\Om,\la)}=\|\cdot\|_{L^\infty(\Om)}$. By invoking \eqref{2.3} and H\"older inequality, we have
\beq \label{2.5}
0\,\leq\, \La_p\,\leq\, \frac{\|\De \xi\|_{L^p(\Om,\la)}}{\|\xi\|_{L^p(\Om,\la)}} \,\leq\, \frac{1}{1-\om(1/p)}\, \frac{\|\De \xi\|_{L^\infty(\Om)}}{\|\xi\|_{L^\infty(\Om)}} \,\leq\, 2\, \frac{\|\De \xi\|_{L^\infty(\Om)}}{\|\xi\|_{L^\infty(\Om)}}.
\eeq
Therefore, by passing to a subsequence, there exists 
\[
\hat{\La}_\infty \in \left[ \, 0\, ,\,\frac{\|\De \xi\|_{L^\infty(\Om)}}{\|\xi\|_{L^\infty(\Om)}}\, \right]
\]
such that $\La_p \larrow \hat{\La}_\infty$, along this subsequence as $p\to \infty$. For technical convenience, we will suppress the subscripts of the subsequences and we will not relabel them. Fix now $k\in \N$, $k\geq 2$. Since $u_p \in \big(W^{2,p}\cap W^{1,p}_0\big)(\Om)$, we have
\[
\ \ \ u_p \in \big(W^{2,k}\cap W^{1,k}_0\big)(\Om), \ \ \ p\geq k.
\]
Since $\p\Om$ is of class $C^{1,1}$, by the Calderon-Zygmund global $L^k$-estimate (see e.g.\ \cite[Lemma 9.17, p.\ 242]{GT}), it follows that there exists a constant $C=C(k,\Om)>0$ such that
\beq \label{2.6}
\phantom{\Big|}   \| u_p \|_{W^{2,k}(\Om)} \,\leq\, C(k,\Om)\, \| \De u_p \|_{L^{k}(\Om,\la)}.
\eeq
By \eqref{2.2},\eqref{2.3} and H\"older inequality, for any $p\geq k$ we have 
\beq \label{2.7}
\phantom{\Big|} \| \De u_p \|_{L^{k}(\Om,\la)} \,\leq\,  \| \De u_p \|_{L^{p}(\Om,\la)} \,=\,\La_p
\eeq
and hence by \eqref{2.5}-\eqref{2.6} we infer that
\beq \label{2.8}
\ \ \ \| u_p \|_{W^{2,k}(\Om)} \, \,\leq\, \frac{2C(k,\Om) \|\De \xi\|_{L^\infty(\Om)}}{\|\xi\|_{L^\infty(\Om)}}, \ \ \ \ p\geq k.
\eeq
By \eqref{2.8}, the sequence $(u_p)_1^\infty$ is bounded in $W^{2,k}(\Om)$ for any $k\in \N$. By passing to a further subsequence if necessary, by Morrey's theorem and a standard weak compactness diagonal argument there exists 
\[
u_\infty \in \bigcap_{1<p<\infty} W^{2,p}(\Om)\cap C^{1}(\overline{\Om})
\]
such that $u_p \larrow u_\infty$ strongly in $C^{1}(\overline{\Om})$ and $\D^2u_p \weak \D^2u_\infty$ weakly in $L^k(\Om,\R^{n\by n})$, as $p\to \infty$ along this subsequence. By the weak lower semicontinuity of the $L^k$ norm, \eqref{2.7} gives for any $k\in \N$ that
\beq \label{eq}
\| \De u_\infty \|_{L^{k}(\Om,\la)} \, \leq\, \liminf_{p\to \infty}\, \| \De u_p \|_{L^{k}(\Om,\la)} \,\leq\, \liminf_{p\to \infty}\,\La_p \, =\, \hat{\La}_\infty
\eeq
and by letting $k\to \infty$, we obtain
\[
\| \De u_\infty \|_{L^\infty(\Om)} \, \leq\, \hat{\La}_\infty.
\]
Further, we have
\[
1 \, =\, \|u_p\|_{L^p(\Om,\la)}\, \leq\, \|u_p\|_{L^\infty(\Om)} \larrow \|u_\infty\|_{L^\infty(\Om)}, \ \text{ as }p\to \infty,
\]
whilst for any $k\in \N$, H\"older's inequality gives
\[
\|u_\infty\|_{L^k(\Om,\la)} \, =\, \lim_{p\to \infty} \, \|u_p\|_{L^k(\Om,\la)}\,\leq\, \limsup_{p\to \infty} \, \|u_p\|_{L^p(\Om,\la)}\,=\, 1.
\]
Therefore, $\|u_\infty\|_{L^\infty(\Om)}=1$. Finally, since $u_\infty \in \big(W^{2,p}\cap W^{1,p}_0\big)(\Om)\set\{0\}$, by \eqref{2.3} and minimality we have
\[
\La_p \,\leq\, \frac{ \|\De u_\infty\|_{L^p(\Om,\la)} }{ \| u_\infty\|_{L^p(\Om,\la)} }
\]
and by letting $p\to \infty$, the above inequality yields $\hat{\La}_\infty \leq \|\De u_\infty\|_{L^\infty(\Om)}$. The lemma ensues. \qed
\ms

The next result shows that the limit $u_\infty$ of $L^p$ minimisers constructed above is an $L^\infty$ minimiser itself.

\begin{lemma} \label{lemma2.4} The function $u_\infty \in \mathcal{W}^{2,\infty}(\Om)$ constructed in Lemma \ref{lemma2.2} solves the minimisation problem
\beq
\label{2.9}
\|\Delta u_\infty\|_{L^\infty(\Omega)} \, =\, \inf \Big\{\|\Delta v\|_{L^\infty(\Om)} \ : \ v \in \mathcal{W}_{\mathrm{H}}^{2,\infty}(\Omega),\, \|v\|_{L^\infty(\Om)} =1 \Big\}
\eeq
where the hinged space $\mathcal{W}_{\mathrm{H}}^{2,\infty}(\Omega)$ is given by \eqref{1.1}. In particular, $\hat{\Lambda}_\infty$ is the value of the infimum in \eqref{2.9}, namely $\Lambda_\infty^{\mathrm{H}} \equiv \Lambda_\infty$. Moreover, we have
\[
\La_\infty \,>\,0.
\]
\end{lemma}

\BPL \ref{lemma2.4}. By \eqref{1.1} we have that $\mathcal{W}_{\mathrm{H}}^{2,\infty}(\Omega) \sub \big({W}^{2,p}\cap W^{1,p}_0\big)(\Omega)$ for all $p\in(1,\infty)$. Hence, by \eqref{2.2}-\eqref{2.3} and minimality, we have
\[
 \La_p \,\leq\,  \|\De v\|_{L^p(\Om,\la)}  , \ \ \ v\in \mathcal{W}_{\mathrm{H}}^{2,\infty}(\Omega), \  \| v\|_{L^p(\Om,\la)}=\,1.
\]
By fixing $v$ and letting $p\to \infty$, by Lemma \ref{lemma2.2} we obtain
\[
\| \De u_\infty \|_{L^\infty(\Om)}\, =\, \hat{\La}_\infty \,\leq\,  \|\De v\|_{L^\infty(\Om)}  , \ \ \ v\in \mathcal{W}_{\mathrm{H}}^{2,\infty}(\Omega), \  \| v\|_{L^\infty(\Om)}=\,1.
\]
By taking infimum over all such $v$, we deduce the equality $\hat{\Lambda}_\infty = \Lambda_\infty$, as claimed. Finally, recall that we already know $0\leq\,  \La_\infty  <\infty$. Suppose for the sake of contradiction that $\La_\infty  =0$. Then, the constraint $ \| u_\infty \|_{L^\infty(\Om)}=1$ contradicts the uniqueness of solutions to the Dirichlet problem for the Laplace equation because  $\De u_\infty =0$ in $\Om$ and $u_\infty =0$ on $\p\Om$. The lemma has been established. \qed
\ms

Next, we prepare towards the construction of the function $f_\infty \in L^1(\Om)\cap BV_{\text{loc}}(\Om)$ and the signed measure $\mu_\infty \in \mathcal{M}(\Om)$ associated with the $\infty$-eigenpair $(u_\infty,\La_\infty)$ which was constructed in Lemmas \ref{lemma2.2}-\ref{lemma2.4} above.

\begin{lemma}\label{lemma2.5} Let $(u_p)_1^\infty$ be the subsequence of the $L^p$ minimisers (satisfying for each $p$ the equalities \eqref{2.2}-\eqref{2.3} and solving the Dirichlet problem \eqref{2.4}) along which the conclusion of Lemmas \ref{lemma2.2}-\ref{lemma2.4} hold. We define the measurable functions $f_p, g_p : \Om \larrow \R$ by 
\begin{align}
f_p\, &:=\, \dfrac{|\De u_{p}|^{{p}-2}\De u_{p}}{(\La_{p})^{p}},
\label{2.10}
\\
g_p\, &:=\, |u_{p}|^{{p}-2}u_{p}.
\label{2.11}
\end{align}
Then, we have
\beq
\label{2.12}
\ \ \De f_p \, =\, g_p \ \ \text{ in }\mD'(\Om),
\eeq
and if $p'=p/(p-1)$, we also have
\begin{align}
\| f_p \|_{L^{p'}(\Om,\la)}\, &=\, \frac{1}{\La_p},
\label{2.13}
\\
\| g_p \|_{L^{p'}(\Om,\la)}\, &=\, 1.
\label{2.14}
\end{align}
\end{lemma}

\BPL \ref{lemma2.5}. The proof is elementary, but we provide it anyway for the sake of completeness. Let $f_p, g_p$ be given by \eqref{2.10}-\eqref{2.11}. We begin by noting that \eqref{2.12} is a consequence of \eqref{2.4} and the definitions. For \eqref{2.13}, by \eqref{2.1}-\eqref{2.3} we have
\[
\begin{split}
\| f_p \|_{L^{p'}(\Om,\la)}\, &=\, \frac{1}{(\La_p)^p} \left(\av_\Om \big| |\De u_p|^{p-2} \De u_p\big|^{\frac{p}{p-1}} \right)^{\frac{p-1}{p}}
\\
&=\, \frac{1}{(\La_p)^p} \left(\av_\Om |\De u_p|^p \right)^{\frac{p-1}{p}}
\\
&=\, \frac{1}{(\La_p)^p} (\La_p)^{p-1}
\\
&=\, \frac{1}{\La_p}
\end{split}
\]
and similarly, in view of \eqref{2.3} we have
\[
\begin{split}
\| g_p \|_{L^{p'}(\Om,\la)}\, &=\,   \left(\av_\Om \big| | u_p|^{p-2} u_p\big|^{\frac{p}{p-1}} \right)^{\frac{p-1}{p}}
\\
&=\,  \left(\av_\Om | u_p|^p \right)^{\frac{p-1}{p}}
\\
&=\, 1.
\end{split}
\]
The lemma ensues. \qed
\ms

\begin{lemma}\label{lemma2.6} In the setting of Lemma \ref{lemma2.5}, there exist a function $f_\infty \in L^1(\Om)\cap BV_{\text{loc}}(\Om)$ and a signed Radon measure $\mu_\infty \in \mathcal{M}(\Om)$ associated with the $\infty$-eigenpair $(u_\infty,\La_\infty)$ such that
\[
\begin{split}
& f_p \larrow f_\infty,  \ \ \ \ \ \, \text{ in }L^q_{\text{loc}}(\Om)\text{ for all }q \in \Big[1,\dfrac{n}{n-1}\Big),
\\
& f_p \weakstar\, f_\infty, \ \ \ \ \ \ \text{ in }BV_{\text{loc}}(\Om),
\\
& g_p \, \mL^n \, \weakstar \, \mu_\infty, \ \ \text{ in }\mathcal{M}(\Om),
\end{split}
\]
along perhaps a further subsequence as $p\to \infty$. Moreover, $f_\infty$ is a distributional solution to the Poisson equation with right hand side $\mu_\infty$:
\[
\ \ \De f_\infty  =\, \mu_\infty \ \  \text{ in }\mD'(\Om).
\]
\end{lemma}

\BPL \ref{lemma2.6}. By Lemmas \ref{lemma2.2} and \ref{lemma2.6}, we have that the sequences $(f_p)_1^\infty$, $(g_p)_1^\infty$ are uniformly bounded in $L^1(\Om)$ and for each $p$ along a subsequence they satisfy
\[
\De f_p \, =\, g_p \ \ \text{ in }\mD'(\Om).
\]
By Lemma \ref{lemma15} and Corollary \ref{corollary16} in our first Appendix, we have that $(f_p)_1^\infty$ is bounded in $L^{n/(n-1)}_{\text{loc}}(\Om)\cap BV_{\text{loc}}(\Om)$ and there exists a limit function $f_\infty$ such that the desired modes of convergence hold true. Moreover, since the absolutely continuous measures $(g_p\, \mL^n)_1^\infty \sub \mathcal{M}(\Om)$ have bounded total variation, there exists a signed Radon measure $\mu_\infty$ such that the desired weak* convergence holds true as well. By passing to the weak* limit in \eqref{2.12} as $p\to \infty$ along an appropriate subsequence, we obtain $\De f_\infty =\mu_\infty$ on $\Om$ in the sense of distributions.

It remains to show that $f_\infty \in L^1(\Om)$. Indeed, fix a compact set $K\sub \Om$ with positive measure. Since $f_p \larrow f_\infty$ as $p\to \infty$ strongly in $L^1_{\text{loc}}(\Om)$ and $(f_p)_1^\infty$ is bounded in $L^1(\Om)$, by \eqref{2.13} and \eqref{2.1} we have
\[
\| f_\infty\|_{L^1(K)} \,=\, \lim_{p\to \infty}\, \| f_p\|_{L^1(K)} \,\leq\, \limsup_{p\to \infty}\| f_p\|_{L^1(\Om)}\,\leq\, \frac{\mL^n(\Om)}{\La_\infty}. 
\]
We conclude by invoking the upper continuity properties of the measure $\| f_\infty\|_{L^1(\cdot)}$ on $\Om$.     \qed
\ms

Now we show the validity of the desired differential inclusion which the $\infty$-eigenpair $(u_\infty,\La_\infty)$ satisfies.

\bl \label{lemma2.7}
Let the quadruple $(u_\infty,\La_\infty, f_\infty, \mu_\infty)$ be as in Lemmas \ref{lemma2.2}-\ref{lemma2.6}. Then, we have
\[
\ \ \De u_\infty(x) \, =\, \La_\infty\, \frac{f_\infty(x)}{|f_\infty(x)|}, \ \text{ a.e. }x\in \Om\set \{f_\infty=0\}.
\]
\el

\BPL \ref{lemma2.7}. By \eqref{2.10}, for any $p$ along the subsequence on which the convergence modes of Lemmas \ref{lemma2.2}-\ref{lemma2.6} hold true, we have
\[
 |\De u_{p}|^{{p}-2}\De u_{p} \, =\, (\La_{p})^{p}f_p, \ \text{ on }\Om.
\]
Since the real function $t\mapsto |t|^{p-2}t$ is strictly increasing and invertible on $\R$ with its inverse given by $s\mapsto s|s|^{-1+\frac{1}{p-1}}$ when $s\neq 0$, we may rewrite the above definition as
\beq \label{2.15}
\ \ \ \De u_p(x)\, =\, (\La_p)^{\frac{1}{1-\frac{1}{p}}} \, \big|f_p(x)\big|^{\frac{1}{p-1}} \frac{f_p(x)}{|f_p(x)|}\, , \ \ \ \text{ a.e. }x\in \Om\set \{f_\infty=0\}.
\eeq
In view of Lemma \ref{lemma2.6}, we may fix $x$ in a subset of $\Om\set \{f_\infty=0\}$ of full measure on which we have $f_p(x) \larrow f_\infty(x)$ as $p\to \infty$ along a subsequence. Then, for $p$ large enough, we have
\[
0\,<\, \frac{|f_\infty(x)|}{2}\,\leq\, |f_p(x)| \,\leq\, 2\, |f_\infty(x)| \,<\, \infty
\]
and hence
\beq \label{2.16}
\ \ \big|f_p(x)\big|^{\frac{1}{p-1}} \larrow 1, \ \text{ as }p\to \infty.
\eeq
By \eqref{2.15}-\eqref{2.16} and Lemmas \ref{lemma2.2}-\ref{lemma2.4} we infer that for a.e.\ $x\in \Om\set\{f_\infty =0\}$ we have
\beq \label{2.17}
\ \ \De u_p(x) \larrow \La_\infty \frac{f_\infty(x)}{|f_\infty(x)|}, \ \text{ as }p\to \infty.
\eeq
By Lemma \ref{lemma2.2} we also have that $(\De u_p)_1^\infty$ is bounded in $L^{q+1}\big(\Om\set\{f_\infty =0\}\big)$ for any $q\in(1,\infty)$ and hence $q$-equi-integrable in $L^{q}\big(\Om\set\{f_\infty =0\}\big)$, because by H\"older's inequality, \eqref{2.2} and \eqref{2.5} we have
\[
\begin{split}
\| \De u_p \|_{L^q(E)} \, &\leq\, \big(\mL^n(E)\big)^{\frac{1}{q(q+1)}}\, \| \De u_p \|_{L^{q+1}(\Om\set\{f_\infty=0\})}
\\
&\leq\, \big(\mL^n(E)\big)^{\frac{1}{q(q+1)}}\, \big(\mL^n(\Om)\big)^{\frac{1}{q+1}}\| \De u_p \|_{L^{p}(\Om,\la)}
\\
&\leq\, \big(\mL^n(E)\big)^{\frac{1}{q(q+1)}}\, \big(\mL^n(\Om)\big)^{\frac{1}{q+1}}\La_p
\\
&\leq\, \left(2  \big(\mL^n(\Om)\big)^{\frac{1}{q+1}}\frac{\|\De \xi\|_{L^\infty(\Om)}}{\|\xi\|_{L^\infty(\Om)}}\right) \big(\mL^n(E)\big)^{\frac{1}{q(q+1)}} \, ,
\end{split}
\]
for any measurable set $E\sub \Om\set \{f_\infty=0\}$ and $p>q+1$. Therefore, by invoking the Vitali convergence theorem (see e.g.\ \cite{FL}), the boundedness of the domain $\Om$ implies
\[
\ \ \ \ \ \ \De u_p \larrow \La_\infty \frac{f_\infty }{|f_\infty |}, \ \text{ in }L^q\big(\Om\set\{f_\infty =0\}\big)
\]
as $p\to \infty$, for any $q\in(1,\infty)$. Since $\De u_p \weak \De u_\infty$ in $L^q\big(\Om\set\{f_\infty =0\}\big)$, uniqueness of weak limits establishes the desired equality and the lemma ensues. \qed
\ms

By combining Lemmas \ref{lemma2.2} and \ref{lemma2.7}, we readily obtain the claimed differential inclusion.

\begin{corollary} \label{corollary2.8}

Let $\mathrm{Sgn} :\R \larrow 2^{\R}$ be the continuous set-valued sign function given by:
\[
\mathrm{Sgn}(x)\,:=\, 
\left\{
\begin{array}{r l}
\{-1\}, & x<0,
\\ 
\phantom{.}[-1,+1], & x=0,
\\
\{+1\}, & x>0.
\end{array}
\right.
\]
Then, the $\infty$-eigenpair $(u_\infty,\La_\infty)$ satisfies
\[
\De u_\infty(x) \, \in\, \La_\infty \mathrm{Sgn}\big(f_\infty(x)\big)\ \  \text{ a.e.\ }x\in\Om.
\]
\end{corollary}

We complete the proof of Theorem \ref{theorem1} by showing that in the case of hinged boundary condition, the differential inclusion reduces to just the Poisson equation with constant right hand side. This result reconciles with the more general findings on (absolute) minimisers of second order functionals in Calculus of Variations in $L^\infty$ in the papers (\cite{MS, S, KP2, KMo}).

\begin{proposition}
\label{proposition3}
Let $\Omega \sub \R^n$ be a bounded domain. Then $u_\infty \in \mathcal{W}^{2,\infty}_{\mathrm{H}}(\Omega)$ is a minimiser for $\Lambda_\infty^{\mathrm{H}}(\Omega)$ if and only if it is a multiple of the solution to
\begin{equation} \label{torsion} 
\left\{
\begin{array}{r l l} -\Delta v \,=\, 1 & \text{in }\Omega, 
\\ 
v \,=\, 0 & \text{on }\partial \Omega. 
\end{array}
\right.
\end{equation}
In particular, $u_\infty$ is strictly positive (or strictly negative) in $\Omega$, and unique up to a nonzero multiplicative constant.
\end{proposition}

Note that for this last part of the proof of the theorem, we do not need any boundary regularity.

\BPP \ref{proposition3}. Let $u_\infty$ be a minimiser realising the infimum in \eqref{1.3}. By a rescaling, we may assume that $\|\Delta u_\infty\|_{L^\infty(\Omega)} = 1$ and by replacing $u_\infty$ by $-u_\infty$, we may assume that 
\[
\| u_\infty\|_{L^\infty(\Om)}\, =\, \underset{\Om}{\ess\sup}\,u_\infty.
\]
Set $g:=-\De u_\infty$ and suppose for the sake of contradiction that $g \not\equiv 1$ on $\Omega$, keeping in mind that $-1\leq g \leq 1$ a.e.\ on $\Om$. To this end, let $v$ be the solution of \eqref{torsion}. We have that
\[
\left\{
\begin{array}{rl}
-\De(v-u_\infty)\, =\, 1-g & \text{ in } \Om,
\ms
\\
v-u_\infty\,=\, 0,\ \ \ \ \, & \text{ on }\p\Om, 
\end{array}
\right.
\]
and $1-g\geq 0$ in $\Om$ with $1-g>0$ on a subset of positive measure. By the strong maximum principle we infer that $u_\infty < v$ in $\Omega$, and therefore
\[
\|u_\infty\|_{L^\infty(\Omega)} \,<\, \|v\|_{L^\infty(\Omega)} 
\]
because the supremum is attained inside $\Om$. This leads to the contradiction to minimality
\[
\frac{\|\De u_\infty\|_{L^\infty(\Om)}}{\|u_\infty\|_{L^\infty(\Om)}} \, >\, \frac{\|\De v\|_{L^\infty(\Om)}}{\|v\|_{L^\infty(\Om)}}.
\]
Therefore, any minimiser $u_\infty$ must satisfy $-\De u_\infty =1$ a.e.\ in $\Om$ up to a scaling. The converse statement can be established by arguing in a completely similar fashion and therefore the conclusion follows.
\qed
\ms

\begin{remark}
We note that Proposition \ref{proposition3} provides existence of a minimiser in the case of hinged boundary conditions, without appealing to the approximation arguments detailed before.
\end{remark}

\BPT \ref{theorem1}. The proof of our first main result is now an immediate consequence of Lemmas \ref{lemma2.2}-\ref{lemma2.7}, Corollary \ref{corollary2.8} and Proposition \ref{proposition3}. \qed
\ms

\ms

\section{The Faber-Krahn inequality for the $\infty$-Bilaplacian and $\infty$-eigenpairs in the case of the ball}

In this section we establish the proof of Theorem \ref{theorem2} in the case of hinged and clamped boundary conditions, whilst we also calculate the eigenvalues and the eigenfunctions in the case that the domain is a Euclidean ball.

\ms

\noi \textbf{The case of hinged boundary conditions.} We begin with the simpler case of hinged boundary conditions. In this section we will be using the symbolisation $\om_n$ for the volume of the unit ball in $\R^n$, whilst $\mathrm{B}_R$ will stand for the open ball in $\R^n$ of radius $R>0$, allowing ourselves the convenient flexibility to mean either centred at the origin, or at any other point. The meaning will be clear from the context and in any case the invariance of the $\infty$-eigenvalue problem under rigid motions will not entail any ramifications.

\begin{proposition}[The Faber-Krahn inequality in the hinged case]
\label{proposition5}
Let $\Omega \sub \R^n$ be a bounded domain with $C^{1,1}$ boundary and let $\mathrm{B}_R$ be a ball with radius
\[
R\,:=\, \left(\frac{\mL^n(\Om)}{\om_n}\right)^{1/n}
\]
namely, such that $\mL^n(\Omega)=\mL^n(\mathrm{B}_R)$. Let $\Lambda_\infty^{\mathrm{H}}(\Om)$ be given by \eqref{1.3}. Then,
\[ 
\Lambda_\infty^{\mathrm{H}}(\Omega) \,\geq \, \Lambda_\infty^{\mathrm{H}}(\mathrm{B}_R),
\]
and equality holds if and only if $\Omega$ coincides with the ball $\mathrm{B}_R$ up to a rigid motion in $\R^n$.
\end{proposition}

\BPP \ref{proposition5}. The proof is a consequence of Talenti's symmetrisation principle \cite[Theorem 1]{talenti}, and of the characterisation of $\infty-$eigenfunctions in Proposition \ref{proposition3}. By a rescaling argument, we may assume without loss of generality that $\mL^n(\Omega)=\mL^n(\mathrm{B}_1)=\omega_n$. Let $u$ be a positive minimiser for $\Lambda_\infty^{\mathrm{H}}(\Omega)$. By \cite[Theorem 1]{talenti}, if $v$ is the solution of the problem 
\[ 
\left\{
\begin{array}{r l l} 
-\Delta v \,=\, 1 & \text{in }\mathrm{B}_1, 
\\ 
v \,=\, 0 & \text{on }\partial \mathrm{B}_1, 
\end{array}
\right.
\]
we obtain that $0 \leq u^* \leq v$ in $\mathrm{B}_1$, where $u^*$ is the Schwarz symmetrisation of $u$. Therefore, we deduce that
\[
\|u\|_{L^\infty(\Omega)} \,= \,\|u^*\|_{L^\infty(\Omega)} \,\leq \,\|v\|_{L^\infty(\Omega)},\]
which implies $\Lambda_\infty^{\mathrm{H}}(\Omega) \geq \Lambda_\infty^{\mathrm{H}}(\mathrm{B}_1)$. By the results of \cite{kesavan2}, it follows that equality holds if and only if $\Omega$ coincides with $\mathrm{B}_1$, up to rigid motions.
\qed
\ms

The next lemma, which is a direct consequence of Proposition \ref{proposition3} of the previous section, completes the picture in the case of hinged boundary conditions.

\begin{corollary}[The $\infty$-eigenpairs in the hinged case]
\label{corollary4}
Let $\mathrm{B}_R$ be the ball of radius $R$ in $\R^n$ centred at the origin. Then every minimiser is a nonzero multiple of the function defined as 
\[ 
u_\infty(x) \,:= \, \frac{1}{2n}(R^2 - |x|^2)
\]
and we also have
\[ 
\Lambda_\infty^{\mathrm{H}}(\mathrm{B}_R) \,= \, \frac{2n}{R^2}.
\]
\end{corollary}

\ms

\noi \textbf{The case of clamped boundary conditions.} We continue with the more complex case of clamped boundary conditions. Let us begin by noting that, if $u \in \mathcal{W}^{2,\infty}_{\mathrm{C}}(\Omega)$, then 
\[
\int_\Omega \Delta u \, =\, 0, 
\]
as a consequence of the Gauss-Green theorem. Nonetheless, the converse is not true in general for a function $u\in \mathcal{W}^{2,\infty}(\Om)$ (satisfying $u=0$ on $\p\Om$), unless $\Omega$ is a ball $\mathrm{B}_R$ and $u$ is radially symmetric. In this case,
\[ 
0 \, =\, \int_{\Omega} \Delta u \,=\, \int_{\partial \Omega} \D u \cdot \nu \, d\mH^{n-1} \,=\, u'(R)\, \mH^{n-1}(\partial \mathrm{B}_R) 
\]
which implies that $u'(R)=0$ and hence indeed $u \in \mathcal{W}^{2,\infty}_{\mathrm{C}}(\Omega)$ as claimed. In the above argument, $\mH^{n-1}$ denotes the $n-1$-Hausdorff measure restricted to $\p\Om$ and $\nu$ the outwards pointing normal vector field on $\p\Om$.

Before proving the Faber-Krahn inequality, we need some technical preparation which is the content of the next lemma.

\begin{lemma} \label{Rdependence}
Let $R \in (0,1]$, and $\mathrm{B}_R \sub \R^n$ be the ball of radius $R$ centred at the origin. Let $f$ be defined on $\mathrm{B}_1$ as
\[ 
f(x) \,:=\, \left\{\begin{array}{r l} 1, & \text{for }|x|\leq 2^{-\frac{1}{n}} ,
\ms
\\ 
-1, & \text{for } 2^{-\frac{1}{n}}<|x|<1 ,
\end{array}
\right.
\]
and let $f_R$ be the restriction of $f$ to $\mathrm{B}_R$. Let $w_R$ be the solution to the problem
\[ 
\left\{
\begin{array}{r l} 
-\Delta w_R \, =\, f_R & \text{in }\mathrm{B}_R, 
\ms
\\ 
w_R \, =\, 0\ \, & \text{on }\partial \mathrm{B}_R.
\end{array}
\right.
\]
Then, when $n = 2$, $w_R$ is given by
\[ 
w_R(x) = \frac{1}{4}(R^2 - |x|^2)
\]
if $R \leq 2^{-\frac{1}{2}}$, and
\[ 
w_R(x) \,= \, 
\left\{
\begin{array}{l l} 
 \dfrac{1}{4} - \dfrac{R^2}{4} + \dfrac{\ln{R}}{2} + \dfrac{\ln{2}}{4} - \dfrac{|x|^2}{4}, & \text{ for }|x|\leq 2^{-\frac{1}{2}} 
\ms\ms
\\ 
 \dfrac{|x|^2}{4} - \dfrac{\ln{|x|}}{2}- \dfrac{R^2}{4} + \dfrac{\ln{R}}{2},   & \text{ for } 2^{-\frac{1}{2}}<|x|<R, 
\end{array}
\right.
\]
otherwise.
If $n \geq 3$, $w_R$ is given by
\[ 
w_R(x) \, = \,\frac{1}{2n}(R^2 - |x|^2)
\]
if $R \leq 2^{-\frac{1}{n}}$, and
\[ 
w_R(x) \,=\, 
\left\{
\begin{array}{l l} 
\displaystyle \frac{2^{-\frac{2}{n}}}{n} - \frac{R^2}{2n} - \frac{R^{2-n}}{n(n-2)} + \frac{2^{1-\frac{2}{n}}}{n(n-2)} - \frac{|x|^2}{2n}, & \text{ for }|x|\leq 2^{-\frac{1}{n}}, 
\ms\ms
\\ 
\displaystyle \frac{|x|^2}{2n} + \frac{|x|^{2-n}}{n(n-2)}- \frac{R^2}{2n} - \frac{R^{2-n}}{n(n-2)},   & \text{ for } 2^{-\frac{1}{n}}<|x|<R, 
\end{array}
\right.
\]
otherwise.

Moreover, in either case $w_R$ has the following properties:
\begin{enumerate}
 \item[(i)] $w_R > 0$ in $\mathrm{B}_R$;
 \item[(ii)] $w_R$ is radially symmetric and radially decreasing;
 \item[(iii)] for $R=1$, $w_1$ belongs to $\mathcal{W}^{2,\infty}_{\mathrm{C}}(\mathrm{B}_1)$;
 \item[(iv)] the function $R \mapsto \|w_R\|_{L^\infty(\Omega)}$, defined on $(0,1]$, attains a strict maximum for $R=1$.
\end{enumerate}
\end{lemma}

The proof of this result is a computation exercise on the use of derivatives in polar coordinates and therefore we refrain from providing the tedious details of it. Now we have:

\begin{proposition}[The Faber-Krahn inequality in the clamped case]
\label{proposition7}
Let $\Omega \sub \R^n$ be a bounded domain with $C^{1,1}$ boundary and let $\mathrm{B}_R$ be a ball with radius
\[
R\,:=\, \left(\frac{\mL^n(\Om)}{\om_n}\right)^{1/n}
\]
namely, such that $\mL^n(\Omega)=\mL^n(\mathrm{B}_R)$. Let $\Lambda_\infty^{\mathrm{C}}(\Om)$ be given by \eqref{1.4}. Then,
\[ 
\Lambda_\infty^{\mathrm{C}}(\Omega) \, \geq \, \Lambda_\infty^{\mathrm{C}}(\mathrm{B}_R),
\]
and equality holds if and only if $\Omega$ coincides with the ball $\mathrm{B}_R$ up to a rigid motion in $\R^n$.
\end{proposition}

\BPP \ref{proposition7}. Without loss of generality, we may assume that $\mL^n(\Omega)=\mL^n(\mathrm{B}_1)=\omega_n$. Let $u$ be a minimiser realising the infimum in \eqref{1.4} for $\Omega$, rescaled in a way that $\|\Delta u\|_{L^\infty(\Omega)} = 1$. By replacing $u$ with $-u$ if necessary, since $u=0$ on $\p\Om$ we may suppose that
\[
\|u\|_{L^\infty(\Om)} \, =\, u(\overline{x}) \, >\, 0
\]
for an interior maximum point $\overline{x} \in \Om$. For convenience we set $f:=-\Delta u$. By the representation formula for solutions of the Poisson equation for the Laplacian (see e.g.\ \cite[Ch.\ 2]{GT}), we have
\[ 
\|u\|_{L^\infty(\Omega)} = u(\overline{x}) = \int_\Omega G(\overline{x},y)f(y)\,dy,
\]
where $G$ is the Green function for $\Omega$. The existence of the latter is guaranteed by the $C^{1,1}$ regularity of the boundary $\p\Om$. By the bathtub principle (\cite[Theorem 1.14]{liebloss}, recalled in our second Appendix), since $\int_\Omega f(y)\,dy = 0$, we have 
\begin{equation} \label{simpadisug} 
u(\overline{x}) \leq \int_\Omega G(\overline{x},y)\widetilde{f}(y)\,dy, 
\end{equation}
with $\widetilde{f} = \chi_{E} - \chi_{\Omega \setminus E}$, where 
\[
E\,:=\,\big\{y \in \Omega \,|\,G(\overline{x},y) > t\big\}
\]
for a suitable $t$ such that $\int_\Omega \widetilde{f}(y)\,dy=0$. Note that we have used the fact that the level sets of $G(\overline{x},\cdot)$ are negligible with respect to the $n$-dimensional Lebesgue measure because $G(\overline{x},\cdot)$ is a harmonic function on $\Omega \setminus \{\overline{x}\}$ (see e.g.\ \cite{HS}). Let now $v$ be the solution of
\[ 
\left\{\begin{array}{r l l}
 -\Delta v \,=\, \widetilde{f} & \text{in }\Omega, 
 \ms
 \\ 
 v \,=\, 0 & \text{on }\partial \Omega.
 \end{array}
 \right.
 \]
Inequality \eqref{simpadisug} reads $0< u(\overline{x})\leq v(\overline{x})$ and therefore 
\[
\|u\|_{L^\infty(\Omega)} \, \leq \, \|v^+\|_{L^\infty(\Omega)}, 
\]
where $v^+$ is the positive part of $v$. Let $\Omega^+$ be the open set $\{v>0\}$ and suppose that $\mL^n(\Omega^+)=\omega_n R^n$. Clearly, we have that $R \in (0,1]$. By Talenti's symmetrisation principle (\cite[Theorem 3.1.1]{kesavan}, recalled in our second Appendix), if $\widetilde{w}_R$ is the solution of the problem
\[ 
\left\{
\begin{array}{r l l} 
-\Delta \widetilde{w}_R \,=\, \widetilde{f}^* & \text{in }\mathrm{B}_R, 
\ms
\\ 
\widetilde{w}_R \,=\, 0 \ \, & \text{on }\partial \mathrm{B}_R,
\end{array}
\right.
\]
then 
\[
\|v^+\|_{L^\infty(\Omega)} \, \leq \, \|\widetilde{w}_R\|_{L^\infty(\Omega)}. 
\]
Let $f_R$ and $w_R$ be the functions defined in Lemma \ref{Rdependence}. By invoking the maximum principle, we obtain that $0 < \widetilde{w}_R \leq w_R$, and thus 
\[
\|v^+\|_{L^\infty(\Omega)} \, \leq \, \| w_R\|_{L^\infty(\Omega)}. 
\]
The last quantity is maximal for $R=1$ by Lemma \ref{Rdependence}. Since $w_1 \in \mathcal{W}^{2,\infty}_{\mathrm{C}}(\mathrm{B}_1)$, we get that $\Lambda_\infty^{\mathrm{C}}(\mathrm{B}_1) \leq \|w_1\|_{L^\infty(\Omega)}^{-1}$ and hence we obtain
\begin{equation} \label{fabkra} \Lambda_\infty^{\mathrm{C}}(\Omega) \geq \Lambda_\infty^{\mathrm{C}}(\mathrm{B}_1).\end{equation}
If equality holds in \eqref{fabkra}, then all the previous inequalities must be equalities. In particular, we have 
\[
\|v^+\|_{L^\infty(\Omega)} \, =\,  \|\widetilde{w}_R\|_{L^\infty(\Omega)}. 
\]
By \cite{kesavan2}, this implies $\Omega^+ = \mathrm{B}_R$ and $v^+ = \widetilde{w}_R$. Moreover, we see that $\widetilde{w}_R \equiv w_R$ and moreover 
$R=1$ by Lemma \ref{Rdependence}. Conclusively, this implies $\Omega = \mathrm{B}_1$ after perhaps a translation. \qed
\ms

By arguing in a fashion similar to that of Proposition \ref{proposition7}, one may further quite easily obtain the following result.

\begin{proposition}
\label{proposition8}
Let $\mathrm{B}_R$ be the ball of radius $R$ in $\R^n$ centred at the origin. Then, the minimiser realising the infimum in \eqref{1.4} in the case of the ball is a positive, radially symmetric function $u$, which satisfies $u(x) = w_1(x/R)$, with $w_1$ as defined in Lemma \ref{Rdependence} for $R=1$. In particular,
\[ 
\Lambda_\infty^{\mathrm{C}}(\mathrm{B}_R) \,= \, 
\left\{
\begin{array}{ll}
\dfrac{\ln{2}}{4R^2}, & \text{ if } n=2,
\ms\ms
\\ 
\dfrac{2^{\frac{2}{n}}(n-2)}{1-2^{\frac{2}{n}-1}}\, \dfrac{1}{R^2} , & \text{ if } n\geq 3.
\end{array}
\right.
\]
\end{proposition}

\begin{remark}
\label{remark9}
It is interesting to notice that $\Lambda_\infty^{\mathrm{C}}(\mathrm{B}_R)$ satisfies
\[ 
\frac{2^{\frac{2}{n}}(n-2)}{1-2^{\frac{2}{n}-1}}\, \frac{1}{R^2} \sim \frac{2n}{R^2} \qquad \text{as }n \to +\infty,
\]
and hence, asymptotically it coincides with the first eigenvalue under hinged boundary conditions. This facts holds true also in the linear case $p=2$. Let $J_\nu$, $I_\nu$ be respectively the Bessel function and the modified Bessel function of the first kind of order $\nu \in \R$. The first eigenvalue of the bilaplacian under Navier boundary conditions is equal to the square of the first eigenvalue of the Laplacian under Dirichlet boundary conditions, and therefore
\[ 
\Lambda_2^H(\mathrm{B}_R) \, =\, \frac{j_{\frac{n}{2}-1}^2}{R^2},
\]
where $j_\nu$ is the first zero of $J_\nu$. On the other hand, the first eigenvalue of the bilaplacian under Dirichlet boundary conditions is given by
\[ 
\Lambda_2^C(\mathrm{B}_R) \, =\,  \frac{k_{\frac{n}{2}-1}^2}{R^2},
\]
where $k_\nu$ is the first zero of $J_\nu I_{\nu +1} + I_\nu J_{\nu+1}$. It can be proven that
\[ 
\frac{\Lambda_2^C(\mathrm{B}_R)}{\Lambda_2^H(\mathrm{B}_R)} \,=\, \frac{k_{\frac{n}{2}-1}^2}{j_{\frac{n}{2}-1}^2} \to 1 \qquad \text{as } n \to +\infty.
\]
The aforementioned results can be found in \cite{ashbaughlaugesen}. At present we do not know whether this property holds true for a general domain $\Omega \sub \R^n$.
\end{remark}

\begin{remark}[The case of $n=1$] 
\label{remark10}
The foregoing reasoning can be applied also to the one-dimensional case of an interval $(-R,R)\sub \R$. When $n=1$, the minimiser is given by $u(x)=w(x/R)$, where $w$ is the piecewise quadratic function
\[
 w(x) \, :=\, 
 \left\{\begin{array}{l l} 
 \displaystyle \frac{1}{4}  - \frac{x^2}{2}, & \text{ for }|x|\leq \frac{1}{2}, 
 \ms
 \\ 
 \displaystyle \frac{x^2}{2} -x +\frac{1}{2}, & \text{ for } \frac{1}{2}<|x|<1 \end{array}
 \right.
 \]
and we also have
\[ 
\Lambda_\infty^{\mathrm{C}}\big( (-R,R) \big) \,= \, \frac{4}{R^2}.
\] 
We note that phenomena of piecewise regularity structure for (absolute) minimisers in $L^\infty$ similar to those arising in the clamped case herein have previously been observed in more general settings of  higher order Calculus of Variations in $L^\infty$ for functionals involving the Laplacian in the papers \cite{MS, S, KP, KMo}.
\end{remark}

\ms

\section{Appendix: Local compactness in $L^1$ and $BV$ of distributional solutions to the Poisson equation with $L^1$ data}

In this appendix we establish an estimate regarding the strong local compactness of $L^1$ distributional solutions to the Poisson equation with $L^1$ right hand side. This result is probably already known in the literature, but since we could not locate a precise reference for this fact, we provide a complete proof for the convenience of the reader.

\begin{lemma} \label{lemma15} Let $\Om \sub \R^n$ be open and bounded and let $u, g \in L^1(\Om)$ be such that
\[
\De u \, =\, g \ \text{ in }\Om,
\]
in the sense of distributions. Then:

\smallskip

\noi \emph{(a)} For any compactly contained $\Om'\Subset \Om$, there is a constant $C=C(\Om,\Om')>0$ such that, 
\[
\big\|u(\cdot +z)-u \big\|_{L^1(\Om')}\, \leq\, C\Big(\|u\|_{L^1(\Om)}\,+\, \|g\|_{L^1(\Om)}\Big)|z|.
\]
for any $0<|z|< \frac{1}{2}\dist(\Om', \p\Om)$.

\smallskip

\noi \emph{(b)} The solution $u \in L^1(\Om)$ belongs to $(L^q_{\text{loc}}\cap BV_{\text{loc}})(\Om)$ for all $q \in [1,1^*]$ and there is a $C=C(q,\Om,\Om')>0$ such that
\[
\big\|[\D u]\big\|(\Om')\,+\, \|u \|_{L^q(\Om')}\, \leq\, C\Big(\|u\|_{L^1(\Om)}\,+\, \|g\|_{L^1(\Om)}\Big),
\]
where $[\D u] \in \mathcal{M}(\Om)$ denotes the measure derivative, $\| [\D u]\|(\cdot)$ is the total variation measure and
\[
1^* =\, \frac{n}{n-1}.
\]
\end{lemma}

As a consequence, we have also:

\begin{corollary} \label{corollary16}
In the setting of Lemma \ref{lemma15}, any sequence $(u_i)_{i=1}^\infty$ of solutions to 
\[
\De u_i \, =\, g_i \ \text{ in }\Om, 
\]
bounded in $L^1(\Om)$ is strongly precompact in $L^q_{\text{loc}}(\Om)$ for $q\in [1, 1^*)$ and weakly* precompact in  $BV_{\text{loc}}(\Om)$ if $(g_i)_{i=1}^\infty$ is also bounded in $L^1(\Om)$. Further, any limit point $u$ such that $u_i \larrow u$ as $i\to \infty$ along a subsequence, solves an equation of the type
\[
\De u \, =\, \mu \ \text{ in }\Om, 
\]
where $\mu \in \mathcal{M}(\Om)$ is a signed Radon measure such that $g_i\,\mL^n \weakstar \mu$ as $i\to \infty$.

\end{corollary}

\BPL \ref{lemma15}. (a) Fix an $\Om'\Subset \Om$. By the local nature of the desired estimate and the properties of the Laplace operator, a mollification argument allows us to assume without harming generality that $u, g \in C^\infty(\overline{\Om})$, $\p\Om$ is piecewise $C^\infty$ and $\De u =g$ classically in $\Om$. (Indeed, if $\De u=g$ on $\Om$, then the standard mollification (as e.g.\ in \cite{E}) yields $\De (u*\eta^\e)=g*\eta^\e$ in an inner $\e$-neighbourhood $\Om_\e$ and we may consider an $\Om''$ such that $\Om'\Subset \Om''\Subset \Om_\e$ whose boundary is piecewise spherical.) Let us also understand $g$ as being extended by zero on $\R^n \set\Om$. By Green's formula (see e.g.\ \cite[Ch.\ 2]{GT}), we decompose $u$ as
\[
u\, =\, h\, +\, \Phi *g, 
\]
where $\De h=0$ in $\Om$ and $\Phi$ is the fundamental solution of the Laplace operator. Then, by setting
\[
R(\Om)\,:=\, 2\,\diam(\Om)\,+\, \dist(\Om',\p\Om)
\]
we estimate
\[
\begin{split}
\big\|(\Phi *g)(\cdot +z)-\Phi *g \big\|_{L^1(\Om')}\, &\leq\, \int_\Om\left|\int_\Om \big[\Phi(x+z-y)-\Phi(x-y)\big]g(y)\, dy\right| dx
\\
&\leq\, |z| \int_\Om\int_\Om \int_0^1\big|\D\Phi(x-y+\la z)\big| |g(y)|\, d\la \,dy\,dx
\\
&\leq\, C\,|z| \int_\Om \int_0^1 \left(\int_{\Om-y+\la z}\frac{dw}{|w|^{n-1}}\right)|g(y)|\, d\la \, dy,
\end{split}
\]
for some $C=C(n)>0$. By using that $x-y+\la z$ lies in the ball $\mathrm{B}_{R(\Om)} \sub \R^n$, when $x,y\in \Om$, $\la \in [0,1]$ and $|z|<\frac{1}{2}\dist(\Om', \p\Om)$, we deduce
\[
\big\|(\Phi *g)(\cdot +z)-\Phi *g \big\|_{L^1(\Om')}\, \leq\, C\,|z| \left( \int_{\mathrm{B}_{R(\Om)}}\frac{dw}{|w|^{n-1}}\right) \|g\|_{L^1(\Om)}.
\]
Further, by using that $h=u-\Phi * g$, Young's inequality implies
\[
\begin{split}
\|h\|_{L^1(\Om)}\, & \leq\, \|u\|_{L^1(\Om)}\,+\, \left\| \int_\Om \Phi(\cdot-y)g(y)\,dy
\right\|_{L^1(\Om)}
\\
 & \leq\, \|u\|_{L^1(\Om)}\,+\, \left\| \int_{\R^n} \Big(\chi_{\Om}(y)\Phi(\cdot-y)\Big) g(y)\,dy
\right\|_{L^1(\Om)}
\\
& \leq\, \|u\|_{L^1(\Om)}\,+\, \mL^n(\Om)\,\|g\|_{L^1(\Om)} \|\Phi\|_{L^1\left(\mathrm{B}_{R(\Om)}\right)}.
\end{split}
\]
Finally, we conclude by putting the pieces together. Let $\Om'_\de$ be the open $\de$-neighbourhood of $\Om'$ with $\de:=\frac{1}{2}\dist(\Om', \p\Om)$. Then, we have $\Om'\Subset \Om'_\de \Subset \Om$ and hence
\[
\begin{split}
\big\|u(\cdot +z)-u \big\|_{L^1(\Om')}\, & \leq\, \big\|(\Phi *g)(\cdot +z)-\Phi *g \big\|_{L^1(\Om')}\,+\, \big\|h(\cdot +z)-h \big\|_{L^1(\Om')}
\\
&\leq\, C\Big(|z|\|g\|_{L^1(\Om)}\,+\, |z|\|\D h\|_{L^\infty(\Om'_\de)}\Big)
\\
&\leq\, C\Big(\|g\|_{L^1(\Om)}\,+\, \|h\|_{L^1(\Om)}\Big)|z|
\end{split}
\]
for some $C=C(n,\Om)>0$, where the last line of the estimate is a consequence of the mean value theorem for harmonic functions and interior derivative estimates (see e.g.\ \cite[Ch.\ 2]{GT}). By the above estimates, we obtain finally
\[
\big\|u(\cdot +z)-u \big\|_{L^1(\Om')}\, \leq\, C\Big(\|u\|_{L^1(\Om)}\,+\, \|g\|_{L^1(\Om)}\Big)|z|.
\]
\smallskip

\noi (b) By the obtained estimate, the difference quotients $(\D^{1,h}u\, \mL^n)_{h\neq 0}$ of $u$ have bounded total variation in the space of Radon measures $\mathcal{M}(\Om')$ and hence by well known arguments
\[
\D^{1,h}u \, \mL^n \, \weakstar \,[\D u] \ \text{ in }\mathcal{M}_{\text{loc}}(\Om,\R^n),
\]
as $h\to 0$. The estimate follows from the weak* lower semi-continuity of the total variation norm and the Sobolev inequality in the BV-space (\cite[Ch.\ 5]{EG}). The lemma ensues.    \qed

\ms

\BPCOR \ref{corollary16}. The result is an immediate consequence of the Fr\'echet-Kolmogorov strong compactness theorem (see e.g.\ \cite[Ch.\ 4]{B}), the Vitali convergence theorem (\cite[Ch.\ 2]{FL}) via an equi-integrability argument similar to that employed in Lemma \ref{lemma2.7} and standard results on the weak* compactness of the spaces of BV functions and Radon measures (\cite[Ch.\ 5]{EG}).   \qed

\ms

\section{Appendix: Some useful results}

In this appendix we collect some useful results which have been utilised earlier in the paper. Some of the results are well-known, and we mention them for the reader's convenience.

\ms

\noi \textbf{Symmetrisations.} Let $\Omega \sub \R^n$ be a bounded domain, and let $f$ be a function in $L^1(\Omega)$. We denote by $\Omega^*$ the ball having the same measure as $\Omega$, and by $f^*$ the Schwarz symmetrisation of $f$, as defined in \cite[Section 1.3]{kesavan}. $f^*$ is a radially symmetric, radially decreasing function defined on $\Omega^*$. It is known that
\[ 
\int_{\Omega^*} u^* = \int_\Omega u,
\]
and also that
\[ 
\|u^*\|_{L^p(\Omega^*)} = \|u\|_{L^p(\Omega)}, \  \  p \in [1,+\infty].
\]
Moreover, if $u \geq 0$ and $u \in W^{1,p}_0(\Omega)$, then $u^* \in W^{1,p}_0(\Omega^*)$, and the \emph{P\'olya-Szeg\"o} inequality holds true:
\[ 
\|\D  u^*\|_{L^p(\Omega^*)} \leq \|\D  u\|_{L^p(\Omega)}, \ \  p \in [1,+\infty].
\]
Symmetrisations are the method of choice in order to prove that the first eigenvalue of the $p$-Laplacian is minimal for the ball, among all domains with fixed volume. Unfortunately this approach does not work for the higher order $L^\infty$ problem we are considering herein because if $u \in W^{2,p}(\Omega)$, it does not in general follow that $u^* \in W^{2,p}(\Omega^*)$. In particular, it does not work even for the pre-limiting case of finite $p$. Nonetheless, the following result of Talenti \cite{talenti} (see also \cite[Theorem 3.1.1]{kesavan}), which turned out to be very useful in the context of higher order problems, is utilised in our proofs in an essential fashion:

\begin{theorem}[Talenti's Symmetrisation Principle] Let $f \in L^2(\Omega)$, and let $u \in W^{1,2}_0(\Omega)$ and $v \in W^{1,2}_0(\Omega^*)$ be the weak solutions of the problems
\[ 
\left\{\begin{array}{r l l} 
-\Delta u \, = \, f & \text{in }\Omega, 
\\ 
u \,=\,  0 & \text{on }\partial \Omega, 
\end{array}
\right. 
\qquad 
\left\{\begin{array}{r l l} 
-\Delta v \,=\, f^* & \text{in }\Omega^*, 
\\ 
v \, =\, 0\ \, & \text{on }\partial \Omega^*. 
\end{array}
\right.
\]
If $u \geq 0$, then $u^* \leq v$ in $\Omega^*$.
\end{theorem}

In particular, by the above result it follows that 
\[
\|u\|_{L^\infty(\Omega)} \,\leq \,\|v\|_{L^\infty(\Omega)}. 
\]
Further, by a result of Kesavan \cite{kesavan2}, equality $\|u\|_{L^\infty(\Omega)} = \|v\|_{L^\infty(\Omega)}$ holds true if and only if $\Omega = \Omega^*$, and $f$ is radially symmetric.

\ms

\noi \textbf{The Bathtub principle.} In our proofs we have also used the following simple measure-theoretic fact, whose proof is a special case of a more general result (see \cite[Theorem 1.14]{liebloss}).

\begin{proposition}
Let $\Omega \sub \R^n$ be a bounded domain and $f \in L^1(\Omega)$ a function such that, for every $t \in \R$, the level set $\{f=t\}$ is a Lebesgue null set. Let $a$, $b$, $\ell \in \R$ be fixed and such that $a \leq \ell\leq b$, and consider the set of functions
\[ 
\mathcal{C}\, :=\, \left\{ g \in L^\infty(\Omega)\,: \,a \leq g \leq b \text{ in }\Omega,\ \av_\Omega g(x)\,dx = \ell \right\}.
\]
Then the supremum in the maximisation problem
\[ \sup_{g \in \mathcal{C}}  \int_\Omega f(x)g(x)\,dx \]
is attained at a function $g \in \mathcal{C}$ of the form 
\[
g \,=\, a \chi_{\{f<t\}} \, +\, b \chi_{\{f \geq t\}}, 
\]
for a suitable $t \in \R$ such that the average of $g$ over $\Om$ is $\ell$.
\end{proposition}

\ms

\ms

\noi \textbf{Acknowledgement.} This work was initiated during a visit of E.P. to the University of Reading in December 2016, partially supported by the Engineering and Physical Sciences Research Council grant EP/N017412/1. Hosting and financing institutions are gratefully acknowledged. N.K.\ would like to thank Craig Evans, Robert Jensen, Roger Moser, Juan Manfredi, Tristan Pryer, Giles Shaw and Zisis Sakellaris for scientific discussions on the topic of higher order $L^\infty$ variational problems. Both authors would like to thank the referee of the paper for the careful and prompt reading of the manuscript, as well as for their constructive comments. 

\ms
\ms

\bibliographystyle{amsplain}

\end{document}